\documentclass[12pt,a4paper]{amsart}
\usepackage{leftidx}
\usepackage{amssymb}
\usepackage{color}
\usepackage{tikz}

\newcommand{\clsp}{\overline{\operatorname{span}}}

\newcommand{\Aut}{\operatorname{Aut}}

\newcommand{\ZZ}{\mathbb{Z}}
\newcommand{\CC}{\mathbb{C}}

\newcommand{\NN}{\mathbb{N}}
\newcommand{\TT}{\mathbb{T}}

\newcommand{\Oo}{\mathcal{O}}
\newcommand{\Kk}{\mathcal{K}}

\newcommand{\Ll}{\mathcal{L}}

\newcommand{\Bb}{\mathcal{B}}

\newcommand{\rinn}[3]{\langle #1,#2 \rangle_{#3}}
\newcommand{\linn}[3]{{}_{#3}\langle #1,#2\rangle}

\newtheorem{theorem}{Theorem}[section]
\newtheorem{corollary}[theorem]{Corollary}

\newtheorem*{lemma*}{Lemma}
\newtheorem{proposition}[theorem]{Proposition}

\theoremstyle{remark}

\newtheorem{example}[theorem]{Example}
\newtheorem{example*}{Example}

\theoremstyle{definition}
\newtheorem{definition}[theorem]{Definition}

\begin{document}

\title[Extensions of Hilbert bimodules]{Extensions of Hilbert bimodules and associated Cuntz-Pimsner algebras}

\author{David Robertson}
\address{Department of Mathematics and Computer Science \\
The University of Southern Denmark \\
Campusvej 55, DK-5230 Odense M, Denmark}
\email{dir@imada.sdu.dk}

\begin{abstract}
We extend the definition of an extension of a right Hilbert module to the setting of Hilbert bimodules and show that an extension of Hilbert bimodules induces an extension of Cuntz-Pimsner algebras. We also study the Cuntz-Pimsner algebra associated to the multiplier bimodule and show that an extension can be realised as a restricted direct-sum bimodule. 
\end{abstract}

\maketitle

\section{Introduction}

A $C^*$-correspondence is a pair $(X,A)$ where $X$ is a right Hilbert $A$-module with a left action of $A$ on $X$. In \cite{Pi}, Pimsner first showed how to associate a $C^*$-algebra to certain $C^*$-correspondences. These $C^*$-algebras are now known as Cuntz-Pimsner algebras. Since then, this definition has been extended to all $C^*$-correspondences by Katsura in his series of papers \cite{Ka1,Ka2,Ka3}.

When a $C^*$-correspondence is simultaneously a left Hilbert module, with a left inner-product, we call it a Hilbert bimodule. In this paper, we extend the definition of an extension of right Hilbert modules due to Baki\'{c} and Gulja\v{s} \cite{BaGu2,BaGu1}, to the setting of Hilbert bimodules. We then show that an extension of Hilbert bimodules induces an extension of Cuntz-Pimsner algebras in a functorial way.

We begin Section \ref{sec:preliminaries} with the definition of  Hilbert bimodules, due to Brown, Mingo and Shen \cite{BrMiSh}. We also recall the definition of the multiplier bimodule $M(X)$, first defined for imprimitivity bimodules in \cite{EcRa}, which is a generalisation of the multiplier algebra for $C^*$-algebras.

In Section \ref{sec:extensions} we define extensions of Hilbert bimodules. This definition is more or less unchanged from the analogous definition for right Hilbert modules given in \cite{BaGu1}, we need only show that this definition also works for bimodules. As an example, we show that the multiplier bimodule $(M(X),M(A))$ is the largest essential extension of a full bimodule $(X,A)$. We also show that under certain assumptions, given an ideal $I \lhd A$ we can find a bimodule $(X_I,I)$ with $(X,A)$ as an extension. Furthermore, this extension will be essential if and only if $I$ is an essential ideal of $A$. We then show that the classification of extensions of right Hilbert modules up to a Busby morphism from \cite{BaGu1} extends to the category of Hilbert bimodules.

We continue in Section \ref{sec:cpalgebras} with the definition of the $C^*$-algebra associated to a Hilbert bimodule. We show that this process is functorial, so morphisms between Hilbert bimodules induce $C^*$-algebra homomorphisms between Cuntz-Pimsner algebras. Furthermore, we show that extensions of Hilbert bimodules induce extensions of the associated Cuntz-Pimsner algebras. We also briefly study the Cuntz-Pimsner algebra $\Oo_{M(X)}$ associated to the multiplier bimodule, and show with an example that in general this is smaller than the multiplier algebra $M(\Oo_X)$.

\section{Preliminaries} \label{sec:preliminaries}

For a $C^*$-algebra $A$, a \emph{right Hilbert $A$-module} is a Banach space $X$ equipped
with a non-degenerate right action of $A$, and an $A$-valued inner-product $\langle \cdot  , \cdot \rangle_X$ satisfying
\begin{enumerate}
 \item $\langle \xi,\eta a \rangle_X = \langle \xi,\eta\rangle_X  a$;
 \item $\langle \eta,\xi\rangle_X  = \langle \xi,\eta\rangle_X ^*$; and
 \item $\langle \xi,\xi\rangle_X  \geq 0$ and $\|\xi\| = \sqrt{\|\langle \xi,\xi\rangle_X \|}$
\end{enumerate}
for all $\xi, \eta \in X$ and $a \in A$.

Likewise, a \emph{left Hilbert $A$-module} is a Banach space $X$ with a non-degenerate left action of $A$ and an $A$-valued inner product $_X\langle \cdot  , \cdot  \rangle$ satisfying analogous relations to those above; i.e.
\begin{enumerate}
 \item $_X\langle a \xi,\eta \rangle = a _X\langle \xi,\eta\rangle  $;
 \item $_X\langle \eta,\xi\rangle  = _X\langle \xi,\eta\rangle ^*$; and
 \item $_X\langle \xi,\xi\rangle  \geq 0$ and $\|\xi\| = \sqrt{\|_X\langle \xi,\xi\rangle \|}$
\end{enumerate}
for all $\xi, \eta \in X$ and $a \in A$.

The following definition of a Hilbert bimodule is originally due to Brown, Mingo and Shen \cite{BrMiSh}.
\begin{definition} \label{definition:bimodule}
We say that the pair $(X,A)$ is a \emph{Hilbert bimodule}, when $X$ is simultaneously a left and right Hilbert $A$-module, and satisfies the relation
\[
 \linn \xi \eta X \zeta = \xi \rinn \eta \zeta X.
\]
\end{definition} 

We say a bimodule $(X,A)$ is \emph{full} if both $\linn X X X \subset A$ and $\rinn X X X \subset A$ are dense. In the literature, a full bimodule is also known as an imprimitivity bimodule.

Given a $C^*$-algebra $A$ and Hilbert bimodules $(X,A)$ and $(Y,A)$, we denote by $\Ll(X,Y)$ the set of all \emph{adjointable} operators from $X$ to $Y$; that is,
linear operators $T:X\to Y$ such that there exists a linear operator $T^*:Y\to X$ called
the \emph{adjoint} of $T$ satisfying
\[
\rinn {T \xi} \eta Y  = \rinn \xi {T^* \eta} X
\]
for all $\xi \in X, \eta \in Y$. If the adjoint $T^*$ exists, it is unique. We will write $\Ll(X)$ for $\Ll(X,X)$.
With the usual operator norm $\|T\|=\sup\{\|Tx\|:\|x\|\leq 1\}$, $\Ll(X)$ is a $C^*$-algebra

For $\xi \in X,\eta\in Y$, define $\theta_{\eta,\xi} \in \Ll(X,Y)$ to be the operator satisfying
\[
 \theta_{\eta,\xi}(\zeta) = \eta \rinn \xi \zeta X.
\]

This is an adjointable operator with $(\theta_{\eta,\xi})^* = \theta_{\xi,\eta}$. We call
\[
 \Kk(X,Y) = \clsp\{\theta_{\eta,\xi}:\xi \in X, \eta\in Y\}
\]
the set of \emph{compact} operators. It is easily seen that $\Kk(X,X) = \Kk(X)$ is a closed two-sided ideal in $\Ll(X)$.

The following simple example will be useful when we are defining the $C^*$-algebras associated to Hilbert bimodules.

\begin{example} \label{example:bimodule}
Let $D$ be a $C^*$-algebra. Then the pair $(D,D)$ is a Hilbert bimodule with left and right actions by given by multiplication and left and right inner products given by
\[
 \linn a b D = ab^* \ \mbox{ and } \ \rinn a b D = a^*b
\]
for $a,b \in D$. It is well-known that there are isomorphisms $\Kk(D) \cong D$ and $\Ll(D) = M(D)$.
\end{example}

\begin{definition}
Let $(X,A)$ and $(Y,B)$ be Hilbert bimodules. A morphism from $(X,A)$ to $(Y,B)$ is a pair of maps $(\psi_X,\psi_A)$ where the map $\psi_X:X \to Y$ is linear and $\psi_A:A \to B$ is a $C^*$-homomorphism such that, for any $\xi, \eta \in X, a \in A$ we have
\begin{enumerate}
\item $\psi_X(\xi a) = \psi_X(\xi) \psi_A(a)$
\item $\psi_X(a \xi) = \psi_A(a) \psi_X(\xi)$
\item $\psi_A(\langle \xi, \eta \rangle_X) = \langle \psi_X(\xi), \psi_X(\eta) \rangle_Y$, and
\item $\psi_A(_X\langle \xi, \eta \rangle) = _Y\langle \psi_X(\xi), \psi_X(\eta) \rangle$.
\end{enumerate}

We say a morphism is \emph{injective} if the map $\psi_A:A \to B$ is injective. In this case, $\psi_X$ is also necessarily injective. For full Hilbert bimodules, \cite[Theorem 2.3]{BaGu3} implies that the converse is also true - that is; $\psi_X$ is injective if and only if $\psi_A$ is injective.

\end{definition}

\begin{definition}
Let $(X,A)$ be a Hilbert bimodule. Define the multiplier $M(X) := \Ll(A,X)$
\end{definition}

The following proposition is already well-known, see \cite{EcKaQuRa} for example. A proof is given here for completeness.

\begin{proposition}
Let $(X,A)$ be a full Hilbert bimodule. Then there exist left and right actions, and left and right inner products such that the pair $(M(X),M(A))$ is also a Hilbert bimodule.
\end{proposition}

\begin{proof}
Firstly, by identifying the multiplier algebra $M(A)$ with $\Ll(A)$ as in Example \ref{example:bimodule}, we can define the right action of $M(A)$ on $M(X)$ simply as composition of operators.

For the left action, define a map $\phi:A \to \Ll(X)$ by
\[
 \phi(a) (\xi) = a \xi.
\]
Since we have assumed that $(X,A)$ is a full Hilbert bimodule, $\phi$ is in fact an isomorphism $A \to \Kk(X)$ and hence extends to an isomorphism $\overline{\phi} :M(A) \to \Ll(X)$. So define the left action by
\[
 (m T)(a) = \overline{\phi}(m)(T(a))
\]
for $m \in M(A), T \in M(X)$ and $a \in A$. 

Define the right inner-product $\langle \cdot  , \cdot \rangle_{M(X)}$ by 
\[
 \rinn S T {M(X)} = S^*T.
\]
For the left inner-product, notice that for any $S,T \in M(X)$ we have $ST^* \in \Ll(X)$. So we may define the left action
\[
 \linn S T {M(X)} := \overline{\phi}^{-1}(ST^*).
\]
We need to see that this satisfies Definition \ref{definition:bimodule}. It has already been proven in \cite{BaGu1} that the right module structure satisfies the appropriate conditions, so we instead concentrate on the left module structure.

Fix $m \in M(A)$and $S,T \in M(X)$.
 Then
\begin{eqnarray*}
\linn {m S} T {M(X)} &=& \overline{\phi}^{-1}((mS)T^*) \\
&=& \overline{\phi}^{-1}(\overline{\phi}(m)ST^*) \\
&=& m\overline{\phi}^{-1} (ST^*) \\
&=& m \overline{\pi}(ST^*).
\end{eqnarray*}
Likewise, for $R, S$ and $T$ in $M(X)$ we have
\begin{eqnarray*}
 \linn R S {M(X)} T &=& \overline{\phi}(\overline{\pi}(RS^*))T \\
 &=& (RS^*)\circ T \\
 &=& R\circ (S^*T) \\
 &=& R \rinn S T {M(X)}
\end{eqnarray*}
as required.
\end{proof}

\section{Extensions of Hilbert bimodules} \label{sec:extensions}

We begin this section with the definition of an extension of a Hilbert bimodule. This is modelled on the definition given for right Hilbert modules in \cite{BaGu1}.

First we need the following concept. For a Hilbert bimodule $(X,A)$, we say an ideal $I$ of $A$ is \emph{invariant with respect to $X$} if $I X = X I$.

\begin{definition}
Let $(X,A)$ be a Hilbert bimodule. We say that a Hilbert bimodule $(Y,B)$ is an \emph{extension} of $(X,A)$ if $B$ is a $C^*$-algebra containing $A$ as an ideal, and there is a morphism of Hilbert bimodules $(\psi_X,\psi_A):(X,A) \to (Y,B)$ such that $\psi_A$ is simply the inclusion of $A$ in $B$ and $A$ is invariant with respect to $Y$. We will call $(Y,B)$ an \emph{essential} extension if $A$ is an essential ideal in $B$.
\end{definition}

\begin{proposition}
Let $(X,A)$ be a Hilbert bimodule, and let $I \lhd A$ be invariant with respect to $X$. Define
\[
 X_I := \clsp\{ \xi a: \xi \in X, a \in I\} \subset X.
\]
Then $(X_I,I)$ is a Hilbert bimodule, with left and right inner-products and left and right actions inherited from $(X,A)$. Furthermore, $(X,A)$ is an extension of $(X_I,I)$, and is essential if and only if $I$ is an essential ideal in $A$.
\end{proposition}

\begin{proof}
We must show that the both left and right actions have range in $X_I$ and both left and right inner products have range in $I$. This is clear for the right action, and for the right inner product this follows easily from $A$-linearity. The left action follows precisely from the assumption that $I$ is invariant with respect to $X$. Finally, for the left inner product, fix $\xi, \eta \in X_I$. Then the invariance of $I$ with respect to $X$ implies that there exists some $i \in I$ and $\zeta \in X$ satisfying $\xi = i \zeta$. So
\[
 \linn \xi \eta X = \linn {i \zeta} \eta X = i \linn \zeta \eta X \in I
\]
as required.

Lastly, given the definition of the inner products and actions, it is clear that the natural inclusion of $(X_I,I)$ inside $(X,A)$ is a morphism of Hilbert bimodules.
\end{proof}

As an easy consequence of this proposition, we get the following.

\begin{corollary}
Let $(X,A)$ be a full Hilbert bimodule. Then the multiplier bimodule $(M(X),M(A))$ is an essential extension of $(X,A)$.
\end{corollary}

\begin{proof}
In  light of the previous proposition and the definition of the bimodule $(M(X),M(A))$, it is enough to show that $A$ is an $M(X)$-invariant ideal of $M(A)$. In fact we show that $M(X)\iota_A(A) = \iota_X(X) = \iota_A(A) M(X)$.
Fix $a \in A$ and $T \in M(X)$. For any $b \in A$ we have
\[
 (T    \iota_A(a))(b) =  T(ab) = T(a)   b = \iota_X(T(a))(b)
\]
so $T    \iota_A(a) = \iota_X(T(a)) \in \iota_X(X)$.
Since $X$ is full, the left action of $A$ on $X$ induces an isomorphism $A \to K(X)$, so a simple calculation shows that
\[
 \iota_A(A) M(X) = \Kk(X) \Ll(A,X) \subset \Kk(A,X)
\]
Then from \cite[Lemma 2.32]{RaWi} we have an isomorphism $\Kk(A,X) \to \iota_X(X)$, so we have shown that we have inclusions
$\iota_A(A) M(X) \subset \iota_X(X)$ and $M(X) \iota_A(A) \subset \iota_X(X)$. For surjectivity, fix an approximate unit $\{e_i\}$ for $A$. Then it is straightforward to show that for any $\xi \in X$
\[
\lim (\iota_A(e_i)    \iota_X(\xi)) =\iota_X(\xi) =  \lim(\iota_X(\xi)    \iota_A(e_i)).
\]
So we have surjectivity as required. 
\end{proof}

The remainder of this section is dedicated to classifying extensions of Hilbert bimodules. This is done in \cite{BaGu2} for right Hilbert modules, and as with the rest of this section we simply show that the same constructions and methods work for Hilbert bimodules.

\begin{proposition}
Given an extension $(Y,B)$ of a full Hilbert bimodule $(X,A)$, there is a unique morphism of Hilbert bimodules 
\[
 (\lambda_Y,\lambda_B):(Y,B) \to (M(X),M(A))
\]
such that $(\lambda_Y,\lambda_B) \circ (\psi_X,\psi_A) = (\iota_X,\iota_A)$. Furthermore, $(\lambda_Y,\lambda_B)$ is injective if and only if $(Y,B)$ is essential.
\end{proposition}

\begin{proof}
We know from \cite[Theorem 1.1]{BaGu1} that there is a morphism of right Hilbert modules
\[
 (\lambda_Y,\lambda_B) : (Y,B) \to (M(X),M(A))
\]
satisfying
\[
 \lambda_Y(\xi)(a) = \psi_X^{-1}(\xi    \psi_A(a)) \\
 \lambda_B(b)(a) = \psi_A^{-1}(b    \psi_A(a)).
\]
It is shown in \cite{BaGu1} that $\lambda_Y(\xi)$ is adjointable $\lambda_Y(\xi)^* \in \Ll(X,A)$ satisfying
\[
 \lambda_Y(\xi)^*(\alpha) = \psi_A^{-1} (\rinn \xi {\psi_X(\alpha)} Y).
\]
We need only check that these maps also preserve the left module structure. For the left inner product, first notice that since we assume that $(X,A)$ is full, we have an isomorphism $\pi:A \to \Kk(X)$, so it is enough to show that for any $\alpha, \beta \in X$ and $\xi, \eta \in Y$ we have
\[
 \pi(\lambda_B(\linn  \xi \eta Y) \linn \alpha \beta X) = \lambda_Y(\xi)\lambda_Y(\eta)^* \theta_{\alpha, \beta}.
\]
Calculating, we see
\begin{eqnarray*}
 \lambda_B(\linn  \xi \eta Y) \linn \alpha \beta ) &=& \psi_X^{-1}(\linn  \xi \eta B    \psi_X(\linn \alpha \beta X)) \\
 &=& \psi_X^{-1}(\linn {\xi    \rinn \eta {\psi_X(\alpha)} Y} {\psi_X(\beta)} Y ) \\
 &=&  \psi_X^{-1}(\linn {\xi    \psi_A(\lambda_Y(\eta)^*(\alpha))} {\psi_X(\beta)} Y ) \\ 
&=& \linn {\lambda_Y(\xi)\lambda_Y(\eta)^*(\alpha)} \beta X)
\end{eqnarray*}
Applying the isomorphism $\pi$ to both sides we get
\[
 \pi(\lambda_B(\linn  \xi \eta Y)\linn \alpha \beta X) = \theta_{{\lambda_Y(\xi)\lambda_Y(\eta)^*(\alpha)},\beta} = \lambda_Y(\xi)\lambda_Y(\eta)^* \theta_{\alpha, \beta}
\]
as required.
For the left action, fix some $\xi \in Y, b \in B$ and $a \in A$. Then
\begin{eqnarray*}
\lambda_Y(b    \xi)(a) &=& \psi_X^{-1}((b   \xi)   \psi_A(a)) \\
&=& \psi_X^{-1}(b   (\xi) \psi_A(a)) \\
&=& \lambda_B(b) \psi_X^{-1}(\xi    \psi_A(a)) \\
&=& \lambda_B(b) \lambda_Y(\xi)(a)
\end{eqnarray*}
as required. So the proof is finished.
\end{proof}

\begin{proposition}
Let $(Y,B)$ be an extension of a Hilbert bimodule $(X,A)$. Then with the structure inherited from $(Y,B)$, the pair of quotient spaces $(Y / X, B/ A)$ is also a Hilbert bimodule. Furthermore, the pair of quotient maps $(q_Y,q_B):(Y,B) \to (Y/X,B/A)$ is a morphism in the category of Hilbert bimodules.
\end{proposition}

\begin{proof}
We begin by defining the left and right actions, and the left and right inner-products. For $q_Y(\xi), q_Y(\eta)$ in $Y/X$, define
\begin{align*}
 & \linn {q_Y(\xi)} {q_Y(\eta)} {Y/X} = q_B(\linn \xi \eta Y) \\
 & \rinn {q_Y(\xi)} {q_Y(\eta)} {Y/X} = q_B(\rinn \xi \eta Y).
\end{align*}
Likewise, for $q_B(b) \in B/A$, define
\begin{align*}
 & q_B(b)    q_Y(\xi) = q_Y(b    \xi) \\
 & q_Y(\xi)    q_B(b) = q_Y(\xi    b).
\end{align*}
Notice that these actions are well-defined since we have assumed that $\psi_A(A)$ is $Y$-invariant, and the non-degenracy of the left and right actions on $X$ implies that $B    \psi_X(X) = \psi_X(X)    B = \psi_X(X)$. The well-definedness of the inner-products is clear from the definition.

Secondly, given $q_Y(\xi), q_Y(\eta)$ and $q_Y(\zeta)$ in $Y/X$, we have
\begin{eqnarray*}
 q_Y(\xi)    \rinn {q_Y(\eta)} {q_Y(\zeta)} {Y/X} &=& q_Y(\xi    \rinn \eta \zeta Y) \\
&=& q_Y(\linn \xi \eta Y    \zeta) \\
&=& \linn {q_Y(\xi)} {q_Y(\eta)} {Y/X}    q_Y(\zeta)
\end{eqnarray*}
as required. So we have a Hilbert bimodule. 
\end{proof}

This construction motivates the following definition.

\begin{definition}
 For a Hilbert bimodule $(X,A)$, we denote the quotient spaces by $Q(X) := M(X)/X$ and $Q(A) := M(A)/A$. We call the Hilbert bimodule $(Q(X),Q(A))$ the \emph{corona bimodule}. When $(X,A)$ is full, for an extension $(Y,B)$ of $(X,A)$ there is a unique morphism $(\delta_Y,\delta_B):(Y/X,B/A) \to (Q(X),Q(A))$ satisfying
\[
 (\delta_Y,\delta_B) \circ (q_Y,q_B) = (q_{M(X)},q_{M(A)}) \circ (\lambda_Y,\lambda_B).
\]
 We call $(\delta_Y,\delta_B)$ the \emph{Busby morphism} associated to the extension $(Y,B)$. It is straightforward to see that this is indeed a morphism.
\end{definition}

Now, suppose we have Hilbert bimodules $(X,A), (Y,B)$ and $(Z,C)$, and morphisms 
\[
 (\psi_X,\psi_A):(X,A) \to (Z,C) \mbox{ and } (\omega_Y,\omega_B):(Y,B) \to (Z,C).
\]
 Then we can form the \emph{restricted direct-sum}
\[
 X \oplus_Z Y := \{(\xi,\eta) \in X \oplus Y: \psi_X(\xi) = \omega_Y(b)\}.
\]
If we also form the pullback $C^*$-algebra
\[
 A \oplus_C B = \{(a,b) \in A \oplus B: \psi_A(a) = \omega_B(b)\}
\]
then the pair $(X\oplus_Z Y, A \oplus_C B)$ with left and right inner-products
\begin{equation}
 \linn {(\alpha,\beta)} {(\xi,\eta)} {X \oplus_Z Y} = (\linn \alpha \xi X , \linn \beta \eta Y ) \\
 \rinn {(\alpha,\beta)} {(\xi,\eta)} {X \oplus_Z Y} = (\rinn \alpha \xi X , \rinn \beta \eta Y )
\end{equation}
and left and right actions
\[
 (a,b)    (\xi, \eta) = (a    \xi, b    \eta) \\
 (\xi, \eta)    (a,b) = (\xi    a, \eta    b)
\]
is a Hilbert bimodule. For a proof of this, see \cite[Proposition 3.2]{RoSz} Furthermore, it is easy to check that the maps 
\[
 (p_X,p_A) : (X,A) \to (X \oplus_Z Y,A \oplus_C B)
\]
and 
\[
 (p_Y,p_B) : (Y,B) \to (X \oplus_Z Y,A \oplus_C B)
\]
which are just the projections onto the first and second coordinates respectively, are morphisms of Hilbert bimodules. Note that these maps are not surjective in general.

\begin{proposition}
Let $(X,A)$ be a full Hilbert bimodule, $(Z,C)$ be a Hilbert bimodule, and $(\delta_Z,\delta_C):(Z,C) \to (Q(X),Q(A))$ be a morphism of Hilbert bimodules. Then there exists an extension $(Y,B)$ of $(X,A)$ whose Busby morphism is $(\delta_Z,\delta_C)$. Furthermore, this extension is essential if and only if $(\delta_Z,\delta_C)$ is injective.
\end{proposition}

\begin{proof}
This result is just the Hilbert bimodule version of \cite[Proposition 3.4]{BaGu2}, and the proof is more or less the same. Namely,
\[
 (Y,B) := (M(X) \oplus_{Q(X)} Z, M(A) \oplus_{Q(A)} C)
\]
satisfies the required properties. The second claim follows directly from \cite[Proposition 3.3]{BaGu2}.
\end{proof}

\section{$C^*$-algebras associated to Hilbert bimodules} \label{sec:cpalgebras}

We begin with the definition of a covariant representation of a Hilbert bimodule.

\begin{definition}
A \emph{covariant representation} of a Hilbert bimodule $(X,A)$ on a $C^*$-algebra $D$ is a morphism $(t_X,t_A) : (X,A) \to (D,D)$, where the pair $(D,D)$ has the structure of a Hilbert bimodule as described in Example \ref{example:bimodule}.
\end{definition}

\begin{definition}
Let $(X,A)$ be a Hilbert-bimodule. The Cuntz-Pimsner algebra $\Oo_X$ associated to $(X,A)$ is defined to be the universal $C^*$-algebra generated by covariant representations. We denote the universal covariant representation of $(X,A)$ by
\[
(T_X,T_A):(X,A) \to (\Oo_X,\Oo_X).
\]
\end{definition}

For the proof of the existence of such a universal algebra, see \cite{Ka1}.

\begin{proposition}
Given a morphism $(\psi_X,\psi_A):(X,A) \to (Y,B)$ there exists a unique $C^*$-homomorphism $\Psi:\Oo_X \to \Oo_Y$ such that the following diagram is commutative.
\begin{center}
\begin{tikzpicture}
 \node (X) at (0,4) {$(X,A)$};
 \node (Y) at (4,4) {$(Y,B)$};
 \node (OX) at (0,0) {$\Oo_X$};
 \node (OY) at (4,0) {$\Oo_Y$};
 \draw[->] (X) to node[auto] {$(\psi_X,\psi_A)$} (Y);
 \draw[->] (OX) to node[auto] {$\Psi$} (OY);
 \draw[->] (X) to node[left] {$(T_X,T_A)$} (OX) ;
 \draw[->] (Y) to node[auto] {$(T_Y,T_B)$} (OY);
\end{tikzpicture}
\end{center}
\end{proposition}

The proof of this proposition can be found in \cite[Proposition 2.9]{RoSz}.

In what follows, we will denote morphisms between Hilbert bimodules as lower-case characters, and the corresponding homomorphism between Cuntz-Pimsner algebras will be denoted by the corresponding upper-case character.

\begin{theorem}  \label{prop:morphismkernel}
Let $(\psi_X,\psi_A):(X,A) \to (Y,B)$ be a morphism of Hilbert-bimodules. Then $(\ker(\psi_X),\ker(\psi_A))$ with left and right actions, and left and right inner products inherited from $(X,A)$ is also a Hilbert bimodule. Furthermore, $\Oo_{ker(\psi_X)} \cong \ker(\Psi)$.
\end{theorem}

\begin{proof} We begin by recalling that the $C^*$-algebra $\Oo_X$ admits a gauge action; that is, there exists a map
\[
 \alpha: \TT \to \Aut(\Oo_X)
\]
such that
\[
 \alpha_z(T_X(\xi)) = z T_X(\xi) \ \mbox{ and } \ \alpha_Z(T_A(a)) = T_A(a)
\]
for any $z \in \TT, \xi \in X$ and $a \in A$. Hence, by \cite[Proposition 10.6]{Ka1}, it is enough to show that the ideal $\ker(\Psi)$ is invariant under the gauge action. This follows easily by noticing that for any $z \in \TT, \alpha_z$ commutes with $\Psi$.
\end{proof}

It is worth noting that the previous theorem is not true in general for arbitrary $C^*$-correspondences, and it is for this reason that this study is restricted to the class of Hilbert bimodules.

We easily get the following corollary.

\begin{corollary} \label{cor:quotient}
Let $(Y,B)$ be an extension of a Hilbert bimodule $(X,A)$. Then the short exact sequence of Hilbert bimodules
\begin{center}
\begin{tikzpicture}
 \node (01) at (-5.5,0) {$0$};
 \node (X) at (-3.5,0) {$(X,A)$};
 \node (Y) at (0,0) {$(Y,B)$};
 \node (Z) at (3.5,0) {$(Y/X,B/A)$};
 \node (02) at (6,0) {$0$};
 \draw[->] (01) to node[auto] {} (X);
 \draw[->] (X) to node[auto] {$(\psi_X,\psi_A)$} (Y);
 \draw[->] (Y) to node[auto] {$(q_Y,q_B)$} (Z);
 \draw[->] (Z) to node[auto] {} (02);
\end{tikzpicture}
\end{center}
induces a short exact sequence of Cuntz-Pimsner algebras
\begin{center}
\begin{tikzpicture}
 \node (01) at (-4,0) {$0$};
 \node (X) at (-2,0) {$\Oo_X$};
 \node (Y) at (0,0) {$\Oo_Y$};
 \node (Z) at (2,0) {$\Oo_{Y/X}$};
 \node (02) at (4,0) {$0$};
 \draw[->] (01) to node[auto] {} (X);
 \draw[->] (X) to node[auto] {$\Psi$} (Y);
 \draw[->] (Y) to node[auto] {$Q$} (Z);
 \draw[->] (Z) to node[auto] {} (02);
\end{tikzpicture}
\end{center}
Consequently, we have an isomorphism $\Oo_{Y/X} \cong \Oo_Y / \Oo_X$.
\end{corollary}

\begin{example}
Let $A$ be a $C^*$-algebra and $\alpha \in \Aut(A)$ be an automorphism. Let $X = A$. Then $(X,A)$ has left and right inner-products and right action as in example \ref{example:bimodule}. The automorphism $\alpha$ define a left action $x \cdot a = \alpha(x)a$, and we have a Hilbert bimodule. It is well-know (see \cite{Pi} for example) that $\Oo_X$ is isomorphic to the crossed-product algebra $A \rtimes_\alpha \ZZ$.

Now suppose we have some ideal $I \lhd A$ with $\alpha(I) = I$ and $AI = IA$. With $W = I$ and considering $(W,I)$ as a Hilbert bimodule as above, it is easily seen that $(X,A)$ is an extension of $(W,I)$. Since $\alpha(I) = I, \alpha$ also defines an automorphim of $A/I$. So Corollary \ref{cor:quotient} tells us that the well-known exact sequence of $C^*$-algebras
\begin{center}
\begin{tikzpicture}
 \node (01) at (-5.5,0) {$0$};
 \node (X) at (-3.5,0) {$I \rtimes_\alpha \ZZ$};
 \node (Y) at (0,0) {$A \rtimes_\alpha \ZZ$};
 \node (Z) at (3.5,0) {$A/I \rtimes_\alpha \ZZ$};
 \node (02) at (6,0) {$0$};
 \draw[->] (01) to node[auto] {} (X);
 \draw[->] (X) to node[auto] {$$} (Y);
 \draw[->] (Y) to node[auto] {$$} (Z);
 \draw[->] (Z) to node[auto] {} (02);
\end{tikzpicture}
\end{center}
comes from the exact sequence of Hilbert bimodules
\begin{center}
\begin{tikzpicture}
 \node (01) at (-5.5,0) {$0$};
 \node (X) at (-3.5,0) {$(W,I)$};
 \node (Y) at (0,0) {$(X,A)$};
 \node (Z) at (3.5,0) {$(X/W,A/I)$};
 \node (02) at (6,0) {$0.$};
 \draw[->] (01) to node[auto] {} (X);
 \draw[->] (X) to node[auto] {} (Y);
 \draw[->] (Y) to node[auto] {} (Z);
 \draw[->] (Z) to node[auto] {} (02);
\end{tikzpicture}
\end{center}
\end{example}

We end this section with some results regarding the Cuntz-Pimsner algebra $\Oo_{M(X)}$ associated to the multiplier bimodule $(M(X),M(A))$.

\begin{proposition} \label{prop:multiplierrep}
Let $(X,A)$ be a full Hilbert bimodule. The there exists an injective covariant representation $(\overline{T_X}, \overline{T_A})$ of $(M(X),M(A))$ on the multiplier algebra $M(\Oo_X)$. Furthermore, $\Oo_{M(X)} \cong C^*(\overline{T_X},\overline{T_A})$.
\end{proposition}

\begin{proof}
We begin by showing that $T_A : A \to \Oo_X$ is non-degenerate; that is $T_A(A)\Oo_X = \Oo_X = \Oo_X T_A(A)$. Given that $\Oo_X$ is generated by the image of $X$ it is enough to show that $T_A(A)T_X(X) = T_X(X) = T_X(X) T_A(A)$. This easily follows from the non-degeneracy of the left and right actions.

Hence, if we think of the pair $(T_X,T_A)$ as a morphism $(X,A) \to (M(\Oo_X),M(\Oo_X))$, \cite[Theorem 1.30]{EcKaQuRa} implies that there is an extended morphism $(\overline{T_X},\overline{T_A}) : (M(X),M(A)) \to (M(\Oo_X),M(\Oo_X))$; i.e. a representation of $(M(X),M(A))$ on $M(\Oo_X)$.

It remains to see that this representation is injective. Notice that $\ker(\overline{T_A}) \lhd M(A)$ satisfies $\ker(\overline{T_A}) \cap A = \{0\}$, and since $A$ is an essential ideal in $M(A)$, we must have $\ker(\overline{T_A}) = \{0\}$. So the representation is injective.

Finally, given the gauge-invariant uniqueness theorem \cite[Theorem 6.4]{Ka3}, to see that we have an isomorphism $\Oo_{M(X)} \cong C^*(\overline{T_X},\overline{T_A})$, it is enough to show that the injective representation $(\overline{T_X},\overline{T_A})$ admits a gauge action. Let $\gamma : \TT \to \Aut(\Oo_X)$ denote the gauge action on $\Oo_X$. Then for fixed $z \in \TT$, it is easily shown that the pair $(\gamma_z \circ T_X, \gamma_z \circ T_A)$ is an injective representation of $(X,A)$ on $\Oo_X$. So it extends to an injective representation $(\overline{\gamma_z \circ T_X}, \overline{\gamma_z \circ T_A})$ of $(M(X),M(A))$. It is easily checked that this representation satisfies
\[
 \overline{\gamma_z \circ T_X}(T) = z T \mbox{ and } \overline{\gamma_z \circ T_A}(m) = m
\]
for all $z \in \TT, T \in M(X)$ and $m \in M(A)$. This shows that the representation admits a gauge action as required.
\end{proof}

As a corollary to this result we have the following.

\begin{corollary}
Let $(X,A)$ be a full Hilbert bimodule and $(Y,B)$ be an extension of $(X,A)$. Then we have an isomorphism
\[
 \Oo_Y \cong \Oo_{M(X)} \oplus_{\Oo_{Q(X)}} \Oo_{Y/X}
\]
where the pull-back is taken along $\Lambda : \Oo_Y \to \Oo_{M(X)}$ and the induced Busby map $\Delta : \Oo_{Y/X} \to \Oo_{Q(X)}$.
\end{corollary}

\begin{proof}
It is well-known that there is an isomorphism
\[
 \Oo_Y \cong M(\Oo_X) \oplus_{Q(\Oo_X)} \Oo_Y / \Oo_X.
\]

The previous proposition implies we have $\Oo_{M(X)} \to M(\Oo_X)$ injective, and Corollary \ref{cor:quotient} says we have an isomorphism $\Oo_Y / \Oo_X \cong \Oo_{Y/X}$. Hence it is enough to show that the image of the Busby map $\Delta$ has image inside $\Oo_{Q(X)}$ and this is easily checked.
\end{proof}

In light of Proposition \ref{prop:multiplierrep}, one may ask whether the injective map $\Oo_{M(X)} \to M(\Oo_X)$ extends to an isomorphism. We end the paper with an example illustrating that this is not the case in general.

\begin{example}
Let $A = X = C_0(\NN)$. For $f \in A$ and $\xi, \eta \in X$ define the left and right actions by
\[
 (f    \xi) (n) = f(n) \xi (n) \mbox{ and } (\xi    f) (n) = \xi(n) f(n+1)
\]
and left and right inner-products
\[
 \linn \xi \eta X (n) = \xi(n) \overline{\eta(n)} \mbox{ and } \rinn \xi \eta X (n) = \overline{\xi(n-1)} \eta(n-1).
\]
\end{example}
For $n \in \NN$, let $\xi_n, f_n$ be the characteristic functions in $X$ and $A$ respectively. Define elements of $\Oo_X$
\begin{align*}
 & E_{ii} = T_A(f_i) \\
 & E_{ij} = T_X(\xi_i) \dots T_X(\xi_{j-1}) \mbox{ for } i < j \\
 & E_{ij} = T_X(\xi_{i-1})^* \dots T_X(\xi_j)^* \mbox{ for } i > j.
\end{align*}
Then $\{E_{ij}: i, j \in \NN\}$ are a system of matrix units generating $\Oo_X$, so we get an isomorphism $\Oo_X \cong \Kk(\ell^2(\NN))$. Let $\{e_i : i \in \NN\}$ be an orthonormal basis for $\ell^2(\NN)$.

Now let $S \in M(X) = \Ll(A,X)$. Then since $S$ is right $A$-linear, we must have
\[
 S(f_n) = \lambda_n \xi_{n-1}
\]
for some $\lambda_n \in \CC$. So we see that $M(X) = C_b(\NN)$ - bounded functions on the natural numbers. We also have $M(A) = C_b(\NN)$. Then we have the extended representation $(\overline{T_X},\overline{T_A})$ of the bimodule $(M(X),M(A))$ on $M(\Kk(\ell^2(\NN))) = \Bb(\ell^2(\NN))$. If $S \in M(X) = C_b(\NN)$,
\[
 \overline{T_X}(S)(e_i) = S(i) e_{i+1}.
\]
is a so-called weighted-shift operator. So Proposition \ref{prop:multiplierrep} implies that $\Oo_{M(X)}$ is isomorphic to the $C^*$-subalgebra of $\Bb(\ell^2(\NN))$ generated by the weighted-shift operators. It is shown in \cite{BuDe} that this is strictly smaller than $\Bb(\ell^2(\NN))$. So in this case, $\Oo_{M(X)} \ncong M(\Oo_X)$.

\bibliographystyle{amsplain} \bibliography{Bibliography}

\providecommand{\bysame}{\leavevmode\hbox to3em{\hrulefill}\thinspace}
\providecommand{\MR}{\relax\ifhmode\unskip\space\fi MR }
\providecommand{\MRhref}[2]{%
  \href{http://www.ams.org/mathscinet-getitem?mr=#1}{#2}
}
\providecommand{\href}[2]{#2}
\begin{thebibliography}{10}

\bibitem{BaGu3}
Damir Baki{\'c} and Boris Gulja{\v{s}}, \emph{On a class of module maps of
  {H}ilbert {$C^*$}-modules}, Math. Commun. \textbf{7} (2002), no.~2, 177--192.
  \MR{1952758 (2003k:46084)}

\bibitem{BaGu2}
\bysame, \emph{Extensions of {H}ilbert {$C^*$}-modules. {II}}, Glas. Mat. Ser.
  III \textbf{38(58)} (2003), no.~2, 341--357. \MR{2052751 (2005f:46112)}

\bibitem{BaGu1}
\bysame, \emph{Extensions of {H}ilbert {$C^*$}-modules}, Houston J. Math.
  \textbf{30} (2004), no.~2, 537--558 (electronic). \MR{2084917 (2005e:46112)}

\bibitem{BrMiSh}
Lawrence~G. Brown, James~A. Mingo, and Nien-Tsu Shen, \emph{Quasi-multipliers
  and embeddings of {H}ilbert {$C^\ast$}-bimodules}, Canad. J. Math.
  \textbf{46} (1994), no.~6, 1150--1174. \MR{1304338 (95k:46091)}

\bibitem{BuDe}
John~W. Bunce and James~A. Deddens, \emph{{$C^{\ast} $}-algebras generated by
  weighted shifts}, Indiana Univ. Math. J. \textbf{23} (1973/74), 257--271.
  \MR{0341108 (49 \#5858)}

\bibitem{EcKaQuRa}
Siegfried Echterhoff, Steve Kaliszewski, John Quigg, and Iain Raeburn, \emph{A
  categorical approach to imprimitivity theorems for {$C^*$}-dynamical
  systems}, Mem. Amer. Math. Soc. \textbf{180} (2006), no.~850, viii+169.
  \MR{2203930 (2007m:46107)}

\bibitem{EcRa}
Siegfried Echterhoff and Iain Raeburn, \emph{Multipliers of imprimitivity
  bimodules and {M}orita equivalence of crossed products}, Math. Scand.
  \textbf{76} (1995), no.~2, 289--309. \MR{1354585 (97h:46093)}

\bibitem{Ka1}
Takeshi Katsura, \emph{A construction of {$C^*$}-algebras from
  {$C^*$}-correspondences}, Advances in quantum dynamics ({S}outh {H}adley,
  {MA}, 2002), Contemp. Math., vol. 335, Amer. Math. Soc., Providence, RI,
  2003, pp.~173--182. \MR{2029622 (2005k:46131)}

\bibitem{Ka2}
\bysame, \emph{On {$C^*$}-algebras associated with {$C^*$}-correspondences}, J.
  Funct. Anal. \textbf{217} (2004), no.~2, 366--401. \MR{2102572 (2005e:46099)}

\bibitem{Ka3}
\bysame, \emph{Ideal structure of {$C^*$}-algebras associated with
  {$C^*$}-correspondences}, Pacific J. Math. \textbf{230} (2007), no.~1,
  107--145. \MR{2413377 (2009b:46118)}

\bibitem{Pi}
Michael~V. Pimsner, \emph{A class of {$C^*$}-algebras generalizing both
  {C}untz-{K}rieger algebras and crossed products by {${\bf Z}$}}, Free
  probability theory ({W}aterloo, {ON}, 1995), Fields Inst. Commun., vol.~12,
  Amer. Math. Soc., Providence, RI, 1997, pp.~189--212. \MR{1426840
  (97k:46069)}

\bibitem{RaWi}
Iain Raeburn and Dana~P. Williams, \emph{Morita equivalence and
  continuous-trace {$C^*$}-algebras}, Mathematical Surveys and Monographs,
  vol.~60, American Mathematical Society, Providence, RI, 1998. \MR{1634408
  (2000c:46108)}

\bibitem{RoSz}
David Robertson and Wojciech Szyma\'{n}ski, \emph{{$C^\ast$}-algebras
  associated to {$C^\ast$}-correspondences and applications to mirror quantum
  spheres}, Illinois J. Math.

\end{thebibliography}

\end{document}